\documentclass[leqno]{amsart}
\usepackage{amsmath,amsfonts,amsthm,amssymb,indentfirst,epic,url}

\newtheorem{theorem}{Theorem}
\newtheorem{lemma}[theorem]{Lemma}
\newtheorem{corollary}[theorem]{Corollary}

\newtheorem{definition}[theorem]{Definition}
\newtheorem{question}[theorem]{Question}

\newcommand{\Q}{\mathbb Q}
\newcommand{\Z}{\mathbb Z}
\newcommand{\eps}{\varepsilon}
\newcommand{\U}{\mathcal{U}}
\renewcommand{\r}{\mathrm}

\begin{document}

\begin{center}
\texttt{Comments, corrections, and related references welcomed!}\\[.5em]
{\TeX}ed \today
\vspace{2em}
\end{center}

\title%
{More Abelian groups with free duals}
\thanks{This preprint is readable online at
\url{http://math.berkeley.edu/~gbergman/papers/}.
}

\keywords{Subgroups of $\Z^\omega$ with whose duals are free abelian
of uncountable rank.
}

\author{George M. Bergman}
\address{University of California\\
Berkeley, CA 94720-3840, USA}
\email{gbergman@math.berkeley.edu}

\begin{abstract}
In answer to a question of A.\,Blass, J.\,Irwin and G.\,Schlitt,
a subgroup $G$ of the additive group $\Z^\omega$ is constructed whose
dual, $\r{Hom}(G,\Z),$ is free abelian of rank $2^{\aleph_0}.$
The question of whether $\Z^\omega$ has subgroups whose duals
are free of still higher rank is discussed, and
some further classes of subgroups of $\Z^\omega$ are noted.
\end{abstract}
\maketitle

\section*{Introduction.}
The additive group $\Z^\omega$
(the countable direct power of $\Z)$ is a nonfree abelian group
$G$ whose rank (maximum number of linearly independent elements)
is the cardinality of the continuum, $2^{\aleph_0};$ but its dual,
$\r{Hom}(G,\Z),$ is known to be free abelian of merely countable rank
\cite{Specker}.
Blass, Irwin and Schlitt \cite{B+I+S}, after examining generally
which subgroups of $\Z^\omega$ have duals
free of countable rank, ask whether a subgroup
of $\Z^\omega$ can have dual free of \emph{uncountable} rank.
Such a subgroup is here constructed.

We begin (\S\ref{S.first_try}) by sketching briefly
an unsuccessful first
try, noting why it fails, and indicating how the difficulty
can be circumvented.
In \S\S\ref{S.framework}-\ref{S.main} we develop the example
that does so.

In \S\ref{S.better?} we ask whether there are
subgroups of $\Z^\omega$ whose duals are free of still larger ranks,
the obvious upper bound being $2^{2^{\aleph_0}}.$
\S\ref{S.examples} notes some classes of subgroups
of $\Z^\omega$ found while thinking about that question.
(For some other unusual subgroups of $\Z^\omega,$
see \cite{ALSC+RG} and papers cited there.)

(A more precise version of the fact that $\r{Hom}(\Z^\omega,\Z)$ is
free of countable rank, which we shall use
below, says that every homomorphism $\varphi:\Z^\omega\to\Z$
factors through the projection to finitely many coordinates.
Regarded as a property of the \emph{codomain} group $\Z$ of
$\varphi,$ this is expressed by saying that $\Z$
is a \emph{slender} abelian group.
There is considerable literature on slender groups, modules,
and other structures; e.g., \cite[\S94]{Fuchs} \cite{RD:bk}.
However, since this note focuses on the \emph{domain} group $G,$
we shall not use that language here.)

I am indebted to John Steel for a helpful observation used
in~\S\ref{S.better?}.

\section{A first attempt.}\label{S.first_try}
I will sketch here my first try at constructing a group with the
desired property.
As mentioned, the more complicated
construction of \S\S\ref{S.framework}-\ref{S.main} will be motivated
as a way of patching up the difficulty with this one.

Clearly, for any countably infinite set $Y,$ the group
$\Z^Y$ is isomorphic to $\Z^\omega;$ so in place of
$\omega$ we will use a countable set that is
more convenient, the set $[0,1]\cap\Q$ of rational numbers
between $0$ and~$1.$

Let $G$ be the group of $\!\Z\!$-valued functions
on $[0,1]\cap\Q$ which are
constant in a neighborhood of each \emph{irrational} $x\in[0,1].$
Such functions can have infinitely many jumps.
(E.g., consider the function which has the value $0$ at $0,$
while for each natural number
$n,$ it has the value $n$ at all rationals in
the subinterval $(2^{-n-1},\,2^{-n}),$ and $-n$ at $2^{-n}$ itself.)

For each irrational $x\in[0,1],$ let $h_x\in\r{Hom}(G,\Z)$
be the map taking every
$g\in G$ to its constant value in the neighborhood of $x.$
If, further, for
each rational $x\in[0,1]$ we define $h_x$ to take $g\in G$
to its value at $x,$ it is not hard to see that the set of
homomorphisms $h_x$ as $x$ ranges over $[0,1]$ are linearly independent.
If they spanned $\r{Hom}(G,\Z),$ we would
have our desired example.

Why might we hope that the $h_x$ would span $\r{Hom}(G,\Z)$?

Consider any $a\in\r{Hom}(G,\Z).$
If we take any decomposition of $[0,1]$ into countably many
(open, closed, or half-open)
subintervals with rational endpoints, then we can define
elements of $G$ independently on different members of this
decomposition.
This observation can be used to construct homomorphisms
$\Z^\omega\to G,$ by starting with an arbitrary $g\in G,$ and
using the $\!n\!$-th entry of an element of $\Z^\omega$ to ``scale''
the output of $g$ on the $\!n\!$-th of our subintervals
(under some fixed enumeration of those subintervals).
Composing such a map $\Z^\omega\to G$ with $a:G\to\Z,$ we get
a homomorphism $\Z^\omega\to\Z.$
But it is known that every homomorphism $\Z^\omega\to\Z$ is a linear
combination of the evaluation maps at \emph{finitely many} coordinates
\cite[Satz~III]{Specker}.
I hoped to use this to prove that the action of $a$ would similarly
be localized at finitely many points $x_1,\dots,x_n\in[0,1],$ and
would in fact be a linear combination of the corresponding
homomorphisms~$h_{x_i}.$

Unfortunately, this fails to be true.
For example, suppose $x$ is an irrational point in $[0,1],$
and $y_0,\ y_1,\ y_2,\ \dots$ a sequence of distinct
points of $[0,1],$ converging to $x.$
I claim
that a homomorphism $G\to\Z$ can be defined by taking each $g\in G$ to
\begin{equation}\begin{minipage}[c]{35pc}\label{d.+-+-}
$(h_{y_0}(g)-h_{y_1}(g))\ +\ (h_{y_2}(g)-h_{y_3}(g))\ +
\ \dots\ +\ (h_{y_{2m}}(g)-h_{y_{2m+1}}(g))\ +\ \dots\ .$
\end{minipage}\end{equation}
Indeed, as $n\to\infty,$ the $y_n$ approach $x,$ so by definition of
$G,$ the sequence of integers $h_{y_n}(g)$ eventually becomes constant.
Hence all but finitely
many of the parenthesized terms of~(\ref{d.+-+-}) are zero, so the
sum~(\ref{d.+-+-}) is defined, and clearly gives a homomorphism.
But it is not a function of the behavior
of $g$ ``at'' finitely many points of $[0,1].$

The trouble with the argument that suggested the opposite
is that the process of prescribing elements of $G$
independently on countably many intervals cannot be carried out
if these intervals converge to an irrational point $x,$
because $x$ itself won't belong to any of the intervals.
If one nevertheless uses such a family of intervals to piece
together functions
on $[0,1]\cap\Q,$ these will not have the property
of ``local constancy'' near $x$ required to give elements of $G.$

The solution we shall take
below is to weaken that local constancy requirement,
and require instead that each element of the group $G$ we will
define be constant on
a ``perforated'' neighborhood of each irrational $x,$ where countably
many ``perforations'' by subintervals with rational endpoints
are allowed, as long as the fraction of the
space they occupy approaches $0$ as we get close enough to~$x.$

This foils the above counterexample, since
we can enclose successive $y_n$ in tiny perforations,
and on these, let $g$ have values independent of its value at $x;$
thus, the rogue homomorphism~(\ref{d.+-+-}) will no longer be defined.
I expected that this kludge, in solving one
problem, would only create worse difficulties.
But it turned out to work, as we shall see in the next three sections.

Between the above sketch and the construction we will present,
there are also a few cosmetic improvements.
Rather than describing elements of $G$ as
functions on $[0,1]\cap\Q,$ we will, for most of our
development, make them functions on $[0,1],$
i.e., give them genuine values at each $x$ in that interval,
and only at the end restrict them to $[0,1]\cap\Q$ to
obtain a solution to the original problem.
With elements of $G$ expressed as functions on $[0,1],$ we will be
able to describe the ``perforated constancy'' condition
as \emph{continuity} in a certain topology
on that interval, finer than the standard topology.
Finally, rather than working, throughout, with the particular structures
of the set $[0,1],$ its standard topology, and that finer topology,
we will posit (in \S\ref{S.framework})
a general situation of a set $X$ with two topologies
related in certain ways, then describe (in \S\ref{S.[0,1]})
our particular choice of $X=[0,1]$ and its two topologies,
and, finally (in \S\ref{S.main}),
prove our freeness result back in the more general situation.

\section{The general context.}\label{S.framework}

Our construction will start with
a set $X$ given with two Hausdorff topologies,
called the \emph{coarse} and the \emph{fine topology},
satisfying the following three conditions.
(Recall that a topological space is \emph{first-countable}
if every point has a countable neighborhood basis.)
\begin{equation}\begin{minipage}[c]{35pc}\label{d.compact}
Under the coarse topology, $X$ is compact and first-countable.
\end{minipage}\end{equation}
\begin{equation}\begin{minipage}[c]{35pc}\label{d.basis}
The coarse topology on $X$ has a basis of open sets
whose members are clopen (closed and open) in the fine topology.
\end{minipage}\end{equation}
\begin{equation}\begin{minipage}[c]{35pc}\label{d.clopen_cover}
For every infinite subset $S$ of $X,$ there exists a
decomposition of $X$ into disjoint subsets clopen
in the fine topology, infinitely many of which contain
members of $S$ in their interiors with respect to the coarse topology.
\end{minipage}\end{equation}

Condition~(\ref{d.basis}) implies that the fine topology contains
the coarse topology, justifying the names.

Conditions~(\ref{d.compact})-(\ref{d.clopen_cover}) would be satisfied
if we took any compact, Hausdorff, first-countable topology
for the coarse topology,
and the discrete topology for the fine topology.
However, our application of this setup will involve
finding a subset $Y\subseteq X$ of comparatively
small cardinality that is dense in the fine topology.
So we will need a fine topology that is not \emph{too} fine.

Note that if we take for $X$ the real unit interval $[0,1],$ and
for the coarse topology the standard topology, and (following the
original idea of the preceding section) define the fine topology to
have for its open sets those sets $U\subseteq[0,1]$ which contain a
neighborhood, with respect to the coarse topology,
of every \emph{irrational} point $x\in U,$ then
conditions~(\ref{d.compact}) and~(\ref{d.basis}) will hold
(the latter because intervals with rational endpoints are clopen in the
fine topology); and the countable
subset $[0,1]\cap\Q$ will be dense in the fine topology.
But~(\ref{d.clopen_cover})
will fail: if we take for $S$ a sequence of points
$y_0,\ y_1,\ y_2,\ \dots$ converging to an irrational point $x,$
then for any covering of $X$ by disjoint sets clopen in the
fine topology, the member of that covering containing $x$
will contain all but finitely many of the $y_n,$ so it is impossible
for infinitely many other members of our covering to meet $S.$
In the next section we shall see that a fine topology
based on ``perforated neighborhoods'' does
satisfy~(\ref{d.clopen_cover}); and in~\S\ref{S.main},
condition~(\ref{d.clopen_cover}) will be the key to constructing
enough composite maps $\Z^\omega\to G\to\Z$ to establish, as was not
true for the $G$ of the preceding section,
that $\r{Hom}(G,\Z)$ is free on the desired basis.

\section{Our two topologies on $[0,1].$}\label{S.[0,1]}

As suggested above, let us take the real unit interval $[0,1]$
for our set $X,$ and
for our coarse topology the usual compact topology.
The description of the fine topology in the next
definition will make precise the ``small perforations'' idea.

In speaking of subintervals of $[0,1],$
we will use the terms ``closed interval'', ``open interval'' and
``half-open interval'' in their conventional senses for real
intervals (of which the first two match
the topological properties of these intervals
under the coarse topology).

\begin{definition}\label{D.fine}
A subset $U\subseteq[0,1]$ will be \emph{open} in the fine
topology if for every irrational $x\in U,$ and every $\eps>0,$
there exists $\delta>0$ such that for every closed interval
$[r,s]\subseteq[0,1]$ with rational endpoints, which contains $x$
and has length $<\delta,$ there exists a finite family
of pairwise disjoint closed subintervals,
also with rational endpoints, which lie in $U\cap [r,s],$
and have total length at least $(1-\eps)(s-r).$
\end{definition}

It is not hard to see that the class of sets so defined
does indeed constitute a topology.
(In verifying closure under pairwise intersection, note that
the $(1-\eps)(s-r)$ condition means
that the part of $[r,s]$ missed by $U$ can be enclosed in
finitely many intervals of length totaling $\leq\eps$ times the
length of $[r,s].$
To establish this property at $x$ for an intersection $U\cap V$ of
sets both satisfying it there, take $\delta$ small enough to get
the same condition in each of $U$ and $V,$
with $\eps/2$ in place if $\eps.)$

(The assumption in Definition~\ref{D.fine} that the intervals named
all have rational endpoints could
be omitted without changing the topology defined, as long as we
keep the restriction that $x$ be irrational, and
specify that $x$ lie strictly between $r$ and $s.$
But the present formulation in terms of intervals with
rational endpoints will be convenient for our purposes.
In statements made below, on the other hand, assumptions
of rational endpoints cannot, in general, be dropped.)

Two easily verified observations:
\begin{equation}\begin{minipage}[c]{35pc}\label{d.intervals_clopen}
Every open, closed, or half-open subinterval of $[0,1]$ with
rational endpoints \textup{(}including the degenerate closed
interval $[r,r]=\{r\}$ for each rational $r\in[0,1])$ is clopen in the
fine topology.
\end{minipage}\end{equation}
\begin{equation}\begin{minipage}[c]{35pc}\label{d.dense_isol}
Under the fine topology, the rational points of $[0,1]$ form
a \emph{dense} set of \emph{isolated} points.
\end{minipage}\end{equation}

Of conditions~(\ref{d.compact})-(\ref{d.clopen_cover})
of the preceding section,~(\ref{d.compact}) is clear for this pair of
topologies, and~(\ref{d.basis}) follows from~(\ref{d.intervals_clopen})
applied to open intervals.
Less trivial is

\begin{lemma}\label{L.clopen_cover}
The pair of topologies on $[0,1]$ described above
satisfies~\textup{(\ref{d.clopen_cover})}.
\end{lemma}

\begin{proof}
Let $S$ be an infinite subset of $[0,1].$
Then it must contain an infinite increasing
or decreasing sequence; assume without loss of
generality that it contains a decreasing sequence
$y_1>y_2>\dots,$ with greatest lower bound $x.$
By dropping enough terms from this sequence, we can assume
that $1>y_1$ and that for every $n,$ we have
\begin{equation}\begin{minipage}[c]{35pc}\label{d.1/2}
$y_{n+1}-x\ \leq\ (y_n-x)/2.$
\end{minipage}\end{equation}

Now let us surround each of our points $y_n$ by an interval
$(r_n,s_n)\subseteq[0,1]$ with rational
endpoints, in such a way that the lengths of these
intervals shrink much faster than the points $y_n$ approach $x:$
\begin{equation}\begin{minipage}[c]{35pc}\label{d./2^n}
$s_n-r_n\ <\ (y_n-x)/2^{n+1}.$
\end{minipage}\end{equation}
In view of~(\ref{d.1/2}) and~(\ref{d./2^n}),
the intervals $(r_n,s_n)$ are disjoint, and
by~(\ref{d.intervals_clopen}) they are clopen in the fine topology.
Now let
\begin{equation}\begin{minipage}[c]{35pc}\label{d.V=}
$V\ =\ [0,1]\,-\,\bigcup_{n\geq 1}\,(r_n,s_n).$
\end{minipage}\end{equation}
Since the $(r_n,s_n)$ are open in the fine topology,
$V$ is closed in that topology.
We claim it is also open.
That the condition of Definition~\ref{D.fine} holds at
all irrational points of $V$ other than,
perhaps, $x,$ is immediate: $V$ contains
a genuine interval about every such point.
Let us show that $V$ also satisfies the required
condition at $x,$ if that is irrational.
(Definition~\ref{D.fine} imposes no such requirement
if $x$ is rational.)

Given $\eps>0,$ we can use~(\ref{d.1/2})
and~(\ref{d./2^n}) to find a $\delta$
such that in every interval $(r,s)$ of length $s-r<\delta$
having rational endpoints and containing $x,$
the lengths of the (infinitely many) intervals $(r_n,s_n)$ meeting
$(r,s)$ sum to less than $\eps/2$ times the length of $(r,s).$
For any such $(r,s),$ let us take an open subinterval
$(r',s')\subseteq(r,s)$
about $x$ with rational endpoints, of length less than $\eps/2$ times
the length of $(r,s).$
This will contain all but finitely many of the $(r_n,s_n),$ hence
the union of $(r',s')$ and those $(r_n,s_n)$ that meet $(r,s)$ will
be a \emph{finite} union of open intervals, of total length less than
$\eps$ times the length of $(r,s);$ hence its complement in $(r,s)$ will
be a finite union of closed intervals of total length
$\geq(1-\eps)(s-r),$ as required.

The sets $(r_n,s_n)$ $(n\geq 1)$ and $V$ thus satisfy the conclusion
of~(\ref{d.clopen_cover}): they are clopen, they form a disjoint
covering of $[0,1],$ and infinitely many of them,
namely the $(r_n,s_n),$ contain members of $S.$
\end{proof}

\section{Describing $G,$ and proving $\r{Hom}(G,\Z)$ free.}\label{S.main}
We now return to the general assumptions of \S\ref{S.framework},
letting $X$ be a set given with any two topologies
satisfying (\ref{d.compact})-(\ref{d.clopen_cover}).
Let $G\subseteq \Z^X$ be the group of all $\!\Z\!$-valued
functions on $X$ that are continuous in the fine topology.
(Throughout this section we understand $\Z$ to have the discrete
topology, so continuity of a function $X\to\Z$ means that the
inverse image of each integer is clopen.)
For each $x\in X,$
\begin{equation}\begin{minipage}[c]{35pc}\label{d.h_x}
Let $h_x:G\to\Z$ be the homomorphism of evaluation at $x.$
\end{minipage}\end{equation}

We shall show below that $\r{Hom}(G,\Z)$ is freely generated
by the elements $h_x.$
First, the easy part.

\begin{lemma}\label{L.indep}
The homomorphisms $h_x$ $(x\in X)$ are linearly independent.
\end{lemma}

\begin{proof}
It clearly suffices to show that for any finite family of
distinct points $x_1,\dots,x_n\in X,$ there exists
$g\in\r{Hom}(G,\Z)$ with $h_{x_1}(g)=1,$ and
$h_{x_m}(g)=0$ for $m=2,\dots,n.$

Since $X$ is Hausdorff in the coarse topology, we can find
a neighborhood of $x_1$ in that topology
containing none of $x_2,\dots,x_n.$
By~(\ref{d.basis}), this will contain a
subneighborhood $U$ of $x_1$ that is clopen in the fine topology.
The characteristic function of $U,$ which is continuous in
the fine topology, gives our desired $g.$
\end{proof}

We now begin the process that will decompose every element
of $\r{Hom}(G,\Z)$ as a linear combination of these maps.
First, a result whose
proof uses nothing specific to $\!\Z\!$-valued functions.
By the \emph{support} of an element $g\in\Z^\omega$ we will mean
$\{x\in X\mid g(x)\neq 0\}.$

\begin{lemma}\label{L.oE}
Let $a$ be a nonzero member of $\r{Hom}(G,\Z).$
Then there exists a point $x\in X$ such that
\begin{equation}\begin{minipage}[c]{35pc}\label{d.everynbd}
For every neighborhood $U$ of $x$ in the coarse topology, there exists
$g\in G$ with support contained in $U$ such that $a(g)\neq 0.$
\end{minipage}\end{equation}
\end{lemma}

\begin{proof}
Suppose the contrary.
Then every $x\in X$ has a neighborhood $U_x$ in the coarse
topology such that
every member of $g$ with support in $U_x$ is in the kernel of $a.$
By~(\ref{d.basis}) we can assume that the sets $U_x$ are also clopen in
the fine topology on $X.$
By compactness of the coarse topology, finitely many of these,
say $U_{x_1},\dots,U_{x_n},$ cover $X.$
Hence the sets
\begin{equation}\begin{minipage}[c]{35pc}\label{d.U-U}
$U_{x_1},\quad U_{x_2}-U_{x_1},\quad \dots\,,\quad
U_{x_n}-(U_{x_1}\cup\dots\cup U_{x_{n-1}})$
\end{minipage}\end{equation}
constitute a covering of $X$ by finitely many pairwise
disjoint sets clopen in the fine topology,
such that every member of $G$ with support in one of them is
in $\ker(a).$

But because these sets are clopen, every
member of $G$ is a sum of members of $G$ with supports in
one or another of them; hence every member of $G$ is
in $\ker(a),$ so $a=0,$ contradicting our hypothesis.
\end{proof}

We now want to prove that for each $a\in G,$ there are
only finitely many $x\in X$ such that~(\ref{d.everynbd}) holds.
The key step will be the first assertion of the next lemma.
(The second assertion will be used later.)
The proof of the lemma will call on the following known result.
\begin{equation}\begin{minipage}[c]{35pc}\label{d.finitely}
\cite[Satz~III]{Specker}
Every homomorphism $\Z^\omega\to\Z$
depends on only finitely many coordinates of its argument,
i.e., can be factored $\Z^\omega\to\Z^n\to\Z,$ where the
first map is given by the projections to some $n$ coordinates,
and the second is an arbitrary homomorphism.
\end{minipage}\end{equation}

\begin{lemma}\label{L.fin_Us}
Suppose $X$ is written as the union of a disjoint family of
sets $U_i$ $(i\in I)$ each clopen in the fine topology, and let
$a\in\r{Hom}(G,\Z).$
Then only finitely many $i\in I$ have the property
\begin{equation}\begin{minipage}[c]{35pc}\label{d.nonzero}
There are elements of $G$ with support in $U_i$ on which $a$
has nonzero value.
\end{minipage}\end{equation}

Assume, further, that $I$ is countable, and let $U$
denote the union of the finitely many $U_i$
satisfying~\textup{(\ref{d.nonzero})}.
Then the value of $a$ at every $g\in G$ is determined
by the restriction of $g$ to $U.$
Equivalently, $a$ has in its kernel all elements of $G$
with support in $X-U.$
\end{lemma}

\begin{proof}
Suppose, in contradiction to the first assertion,
that there are infinitely many $i$ satisfying~(\ref{d.nonzero}).
Then we can write $I$ as the union of a countably infinite
family of pairwise disjoint
subsets $I_n$ $(n\in\omega)$ each containing at least one $i$
that satisfies~(\ref{d.nonzero}).
Letting $V_n=\bigcup_{i\in I_n} U_i$ for each $n\in\omega,$ it follows
that for each $n\in\omega$ there exists a $g_n\in G$ with
support in $V_n$ such that $a(g_n)\neq 0.$

Since the $g_n$ have disjoint supports, we see
that for every $f\in\Z^\omega,$ the expression
$f'=\sum_{n\in\omega} f(n)\,g_n$ makes sense; and as those
supports are clopen in the fine topology, $f'$ is again
continuous in that topology, i.e., belongs to $G.$
Clearly the map $f\mapsto f'$ is a homomorphism $\Z^\omega\to G.$
Composing it with the given homomorphism
$a:G\to\Z$ we get a homomorphism
$\Z^\omega\to\Z$ which for each $n$ takes the element
$e_n\in\Z^\omega$ having a $1$ in the $\!n\!$-th position
and $0$ everywhere else to $a(g_n)\neq 0.$
But this contradicts~(\ref{d.finitely}), proving the finiteness
of the set of $i\in I$ such that~(\ref{d.nonzero}) holds.

Now assume, as in the last paragraph of the lemma, that $I$
is countable; without loss of generality we shall take $I=\omega.$
Let $U$ be the union of the finitely many $U_i$ for
which~(\ref{d.nonzero}) holds.
Because $U$ and $X-U$ are clopen in the fine topology, every
element of $g$ is the sum of an element with support
in $U$ and an element with support in $X-U.$
From this, the equivalence of the last two sentences
of the lemma is clear; we shall prove the last of those sentences.

Let $g\in G$ be an element
with support in $X-U,$ and now for each $n\in\omega$ let
$g_n\in G$ be the function which agrees with $g$ on $U_n$
and is zero elsewhere.
Again, we map $\Z^\omega$ to $\Z$ by
$f\mapsto a(\sum_{n\in\omega} f(n)\,g_n).$
This function is zero on every $e_n$ $(n\in\omega):$ on those with
$U_n\subseteq U$ because $g$ has support in $X-U,$ and on
the others because the corresponding
sets $U_n$ do not satisfy~(\ref{d.nonzero}).
But from~(\ref{d.finitely}) we can see that
a homomorphism $\Z^\omega\to\Z$ which is zero on all
the $e_n$ is zero; hence the above map is zero, hence
$a(g),$ which is the value of that map on the constant
function $1,$ is $0,$ as claimed.
\end{proof}

(We shall see in Corollary~\ref{C.main+} that the countability
condition in the second paragraph of the above lemma can be dropped;
but the above version suffices for the purposes of this section.)

We deduce

\begin{lemma}\label{L.finite}
For any $a\in\r{Hom}(G,\Z)$ there are only finitely many
$x\in X$ such that~\textup{(\ref{d.everynbd})} holds.
\end{lemma}

\begin{proof}
If the set $S$ of such points were infinite,
then by~(\ref{d.clopen_cover}) we could find a covering of $X$
by disjoint subsets $U_i$ clopen in the fine topology,
infinitely many of which contained a point of $S$
in their interiors with respect to the coarse topology.
It follows from our choice of $S$ that for each of the
latter sets, we could
find an element $g_i\in G$ with support in $U_i$ such
that $a(g_i)\neq 0,$ contradicting the first assertion
of the preceding lemma.
\end{proof}

From Lemmas~\ref{L.oE} and~\ref{L.finite}, we can now get

\begin{corollary}\label{C.finite_decomp}
Every $a\in\r{Hom}(G,\Z)$ can be written
\begin{equation}\begin{minipage}[c]{35pc}\label{d.a=}
$a\ =\ a_1+\dots+a_n$ $(n\geq 0;\ a_1,\dots,a_n\in\r{Hom}(G,\Z))$
\end{minipage}\end{equation}
where for each $a_m$ there is an $x_m\in X$ which is the unique
point such that~\textup{(\ref{d.everynbd})} holds with
$a_m$ and $x_m$ in the roles of $a$ and $x;$ equivalently, such that
\begin{equation}\begin{minipage}[c]{35pc}\label{d.ann_far}
$a_m$ is nonzero, but annihilates all elements of $G$ whose
supports do not have $x_m$ in their closure in the
coarse topology.
\end{minipage}\end{equation}
\end{corollary}

\begin{proof}
Let $x_1,\dots,x_n$ $(n\geq 0)$ be the points described by
Lemma~\ref{L.finite}.
Using~(\ref{d.basis}), we can
get a covering of $X$ by disjoint sets $U_1,\dots,U_n$ which are clopen
in the fine topology, and such that each $U_m$ is
a neighborhood of $x_m$ in the coarse topology
(cf.\ the method used to construct~(\ref{d.U-U})).
If we define $a_m:G\to\Z$ for $1\leq m\leq n$ to be the operation that
first multiplies $g\in G$ by the characteristic function of $U_m,$
then applies $a$ to the result, we immediately
have~(\ref{d.a=}), and it is not hard to verify that
$x_m$ is the unique point for which $a_m$ satisfies~(\ref{d.everynbd}).

It remains to show that this is equivalent to~(\ref{d.ann_far}).
One direction, that~(\ref{d.ann_far}) implies
that $x_m$ is the unique point which, together with the homomorphism
$a_m,$ satisfies~(\ref{d.everynbd}), is easily checked.

Conversely, let us assume that condition and deduce~(\ref{d.ann_far}).
The assumption that there exists $x_m$ which, with $a_m,$
satisfies~(\ref{d.everynbd}) clearly implies that $a_m\neq 0.$
Now let $g\in G$ be any element whose support does not have
$x_m$ in its closure under the coarse topology.
Thus, $x_m$ has a neighborhood $V$ in that
topology disjoint from the support of $g.$
By~(\ref{d.basis}), $V$ has a subset $W$
which is again a neighborhood of $x_m$ in the course topology,
and which is clopen in the fine topology.
The operation of multiplying by the characteristic function
of $X-W$ and then applying $a_m$ will be a member
of $\r{Hom}(G,\Z)$ having no point $x$ satisfying the
analog of~(\ref{d.everynbd}), since the characteristic
function of $W$ shows that $x_m$ is not such a point,
and the fact that $a_m$ has no such point other than $x_m$
shows that the function constructed from it can't either.
Hence by Lemma~\ref{L.oE}, this new function is the zero map.
Finally, because $g$ has support in $X-W,$ this map, which we have
shown to be zero, agrees with $a_m$ at $g;$ so $a_m(g)=0,$ as claimed.
\end{proof}

In view of the above results, we will have what we
have been aiming for, once we prove

\begin{lemma}\label{L.h_x}
Suppose $a\in\r{Hom}(G,\Z)$ is a homomorphism for which there
exists a unique $x\in X$ satisfying~\textup{(\ref{d.everynbd})}.
Then $a$ is an integer multiple of $h_x.$
\end{lemma}

\begin{proof}
The asserted conclusion is clearly equivalent to the statement
that $a$ can be factored $G\to\Z\to\Z,$ where the first map
is $h_x,$ and the second is a homomorphism.
Because $h_x$ maps surjectively to $\Z,$ this
is in turn equivalent to the statement that $a$ annihilates
the kernel of $h_x.$
That is what we shall now prove.

Let $g\in\ker(h_x).$
Since $g$ is continuous in the fine topology on $X,$ and $\Z$ is
discrete, the statement that $g$ belongs to $\ker(h_x)$ says that
\begin{equation}\begin{minipage}[c]{35pc}\label{d.0_on_U}
$g$ is zero on some subset $U$ of $X$ which contains $x$ and is clopen
in the fine topology.
\end{minipage}\end{equation}

On the other hand, note
that by~(\ref{d.basis}), $x$ has a neighborhood basis in the
coarse topology consisting of subsets clopen in the fine topology.
Since the coarse topology is first-countable, that neighborhood basis
can be assumed countable.
By taking successive intersections of its terms, we can assume
it is decreasing, and begins with the whole space:
\begin{equation}\begin{minipage}[c]{35pc}\label{d.chain}
$X=V_0\ \supseteq\ V_1\ \supseteq\ \dots\ \supseteq\ V_n
\ \supseteq\ \dots\,.$
\end{minipage}\end{equation}
Since $X$ is Hausdorff, $\bigcap V_n=\{x\}.$
Hence the sets
\begin{equation}\begin{minipage}[c]{35pc}\label{d.V-V}
$V_0-V_1,\quad V_1-V_2,\quad\dots\,,\quad V_n-V_{n+1},\quad\dots$
\end{minipage}\end{equation}
which are also clopen in the fine topology,
give a disjoint covering of $X-\{x\}.$

Now for the tricky step.
Using the $U$ of~(\ref{d.0_on_U}), we obtain from~(\ref{d.V-V}) the sets
\begin{equation}\begin{minipage}[c]{35pc}\label{d.V-V-U}
$U;\quad V_0-V_1-U,\quad V_1-V_2-U,\quad\dots\,,\quad
V_n-V_{n+1}-U,\quad\dots$
\end{minipage}\end{equation}
which are again disjoint, and clopen in the fine topology,
and which now cover all of $X.$
Hence in view of the final paragraph of Lemma~\ref{L.fin_Us},
$a(g)$ is determined by the values of $a$
on the projections of $g$ to these sets.
By the choice of $U$ in~(\ref{d.0_on_U}),
$g$ has zero projection to $U.$
On the other hand, its
projection on each set $V_n-V_{n+1}-U$ is (as a function
on $X)$ zero on a
neighborhood of $x$ in the coarse topology, namely $V_{n+1}.$
But since $x$ is the unique point which,
with $a$ satisfies~(\ref{d.everynbd}), the equivalence
of the two conclusions of the preceding corollary show
that our $g$ is annihilated by $a,$ as required.
\end{proof}

The above results give

\begin{theorem}\label{T.main}
Suppose $X$ is a set given with two topologies,
called ``the coarse topology'' and ``the fine topology'',
satisfying~\textup{(\ref{d.compact})-(\ref{d.clopen_cover})},
and that $G$ is the subgroup of $\Z^X$ consisting of all elements
continuous with respect to the fine topology on $X$ and
the discrete topology on $\Z.$

Then $\r{Hom}(G,\Z)$ is the free abelian group on the
evaluation maps $h_x$ of~\textup{(\ref{d.h_x})}.\qed
\end{theorem}

\begin{corollary}\label{C.main}
Letting $[0,1]$ denote the real unit interval, and $G$ the
group of functions $[0,1]\to\Z$ continuous in the
topology of Definition~\ref{D.fine}, the group $\r{Hom}(G,\Z)$
is free abelian on the generators $h_x$ $(x\in[0,1]).$

Hence, the isomorphic group $G_0$ of functions $[0,1]\cap\Q\to\Z$
obtained from $G$ by restriction to $[0,1]\cap\Q$ has the same dual.

Hence $\Z^\omega$ has a subgroup whose dual is free abelian of continuum
rank.
\end{corollary}

\begin{proof}
The first assertion follows from the above theorem, since we
showed in the preceding section that the standard topology and
the topology of Definition~\ref{D.fine}
satisfy~(\ref{d.compact})-(\ref{d.clopen_cover}).
The one-one-ness of the restriction map $G\to\Z^{[0,1]\cap\Q},$
on which the second assertion
then hangs, follows from the density statement of~(\ref{d.dense_isol}).
Finally, any bijection between $[0,1]\cap\Q$ and $\omega$ induces an
isomorphism between $\Z^{[0,1]\cap\Q}$ and $\Z^\omega,$ so our subgroup
$G_0\subseteq\Z^{[0,1]\cap\Q}$
leads to an isomorphic subgroup of $\Z^\omega.$
\end{proof}

\section{Can we do better?}\label{S.better?}
The group $\Z^\omega$ has continuum cardinality, $c=2^{\aleph_0};$
hence, although it has relatively few homomorphisms to $\Z$
(only countably many), there is no reason why it should not
have subgroups admitting as many as $2^c$ such homomorphisms.

And, in fact, it does.
It is known that the subgroup $B$ of
\emph{bounded} functions $\omega\to\Z$ is free
(proved in \cite{Specker} assuming the continuum hypothesis,
then in \cite{Noebeling} and \cite{ZB} without that assumption).
Having continuum cardinality and being free, $B$ must be free
of continuum rank; so its dual $\r{Hom}(B,\Z)$ can be
identified with $\Z^c,$ and so has cardinality $2^c.$
Hence that dual group has rank (maximum number of linearly
independent elements) also $2^c,$ in other words,
it contains a subgroup free of that rank.

Perhaps unexpectedly, one can characterize an explicit family
of $2^c$ linearly independent elements in $\r{Hom}(B,\Z).$
The set $\beta(\omega)$ of all ultrafilters on $\omega$
has cardinality $2^c$ \cite[Corollary~7.4]{C+N},
and since each $g\in B$ assumes only finitely
many values in $\Z,$ every ultrafilter $\U\in\beta(\omega)$
gives a way of associating to each such $g$ one of those values,
the value such that the set on which it is assumed belongs to $\U.$
For each $\U\in\beta(\omega)$
this gives a homomorphism $h_\U:B\to\Z,$ and it is not hard to
check that these homomorphisms are linearly independent.

But the free group generated by the maps
$h_\U$ is not the whole of $\r{Hom}(B,\Z)\cong
\r{Hom}(\bigoplus_c\Z,\,\Z)\cong\Z^c,$ since
the latter, having subgroups isomorphic to $\Z^\omega,$ is non-free.
So we may ask

\begin{question}\label{Q.2^c}
Is there a subgroup $G\subseteq\Z^\omega$ whose dual is free
of rank $2^c;$ or, at least, free of rank $>c$?
\end{question}

If we had a candidate subgroup $G$ for the above property, somewhat
along the lines of the construction of
the preceding sections, a possible difficulty
with applying the results of those sections to it is that the
coarse topology on $X$ might no longer be first-countable,
as assumed in~(\ref{d.compact}).
Let us show, therefore, that the general results of that section
remain true if hypotheses~(\ref{d.compact})
and~(\ref{d.basis}) are replaced by
\begin{equation}\begin{minipage}[c]{35pc}\label{d.compact+}
Under the coarse topology, $X$ is compact; and the cardinality
of the set $X$ is less than every countably measurable cardinal
(if any such cardinals exist).
\end{minipage}\end{equation}
\begin{equation}\begin{minipage}[c]{35pc}\label{d.basis+}
For every $x\in X,$ there exists a family $(W_i)_{i\in I_x}$ of pairwise
disjoint sets each clopen in the fine topology, such that\\[.2em]
(i) \ $X-\{x\}\ =\ \bigcup_{i\in I_x} W_i,$\\
(ii) \ $x$ has a basis of open neighborhoods $U$ in the
coarse topology, each of which is
clopen in the fine topology, and has the form
$U=\{x\}\cup\bigcup_{i\in J} W_i$ for some subset $J\subseteq I_x.$
\end{minipage}\end{equation}

For the concept of a countably measurable cardinal (in many
works simply called a measurable cardinal) see \cite{Ch+Keis}.
Condition~(\ref{d.compact+}) is indeed a consequence
of~(\ref{d.compact}), since by~\cite{Arhangelskii}, every
compact Hausdorff first-countable topological space has
cardinality $\leq c,$ which is far less than any countably
measurable cardinal.
Condition~(\ref{d.basis+}) by itself
is stronger than~(\ref{d.basis}), but the
construction of~(\ref{d.V-V}) in the proof of Lemma~\ref{L.h_x} shows
that~(\ref{d.compact})-(\ref{d.clopen_cover})
together imply~(\ref{d.basis+}), so the conjunction
of~(\ref{d.compact+}), (\ref{d.basis+}) and~(\ref{d.clopen_cover})
is implied by that of~(\ref{d.compact}), (\ref{d.basis})
and~(\ref{d.clopen_cover}).
Moreover, this new set of conditions does not imply
first-countability of the coarse topology
(as can be verified using examples where the fine
topology is discrete); so it is strictly weaker than the old one.

We can now prove

\begin{corollary}[to the proofs of \S\ref{S.main}]\label{C.main+}
The general results of \S\ref{S.main} \textup{(}the results
through Theorem~\ref{T.main}\textup{)} remain true
if the hypotheses~\textup{(\ref{d.compact})-(\ref{d.clopen_cover})}
are weakened to~\textup{(\ref{d.compact+})},
\textup{(\ref{d.basis+})} and~\textup{(\ref{d.clopen_cover})}.

Moreover \textup{(}under these weakened
hypotheses, and hence under the original hypotheses\textup{)}
one can delete the countability assumption from
the second paragraph of Lemma~\ref{L.fin_Us}.
\end{corollary}

\begin{proof}[Sketch of proof]
The first-countability condition of~(\ref{d.compact}) was not
used before the proof of Lemma~\ref{L.h_x}, so
the results proved up to that point remain true
under our weakened hypotheses, and the only thing we
have to prove regarding those results is
that in the second paragraph of Lemma~\ref{L.fin_Us}, we can
replace the condition that the index set $I$ be countable by the
assumption from~(\ref{d.compact+})
that $X$ have cardinality less than all
countably measurable cardinals.

The latter assumption certainly implies that (after dropping
from $I$ any $i$ such that $U_i$ is empty) the cardinality of $I$ is
likewise less than all countably measurable cardinals.
Our earlier proof of the desired statement
then goes over if, where we previously called on Specker's
result~(\ref{d.finitely})
that $\r{Hom}(\Z^\omega,\Z)$ is free on the projections to the
individual coordinates, we now call on the stronger
known result that the same is
true of $\r{Hom}(\Z^Y,\Z)$ for any set $Y$ having cardinality less than
all countably measurable cardinals
(cf.\ \cite[Theorem~94.4]{Fuchs}, \cite{E+L}).

Moving on to the proof of Lemma~\ref{L.h_x},
the first-countability condition was finally used there to
construct the family of sets $V_n-V_{n+1}.$
We claim that in our present
context, the $W_i$ of~(\ref{d.basis+}) can serve in the same role.
Without loss of generality, let us assume all those $W_i$ nonempty.

In view of the strengthened final
statement of Lemma~\ref{L.fin_Us} noted
above, we do not need to assume
that there are only countably many $W_i;$
the hypothesis on the cardinality of $X$ implies the very weak
bound we need.
The proof of Lemma~\ref{L.h_x} also used
the fact that each $V_n-V_{n+1}$
was disjoint from some neighborhood of $x$ under the coarse
topology, namely, $V_{n+1}.$
To get the same conclusion for a given $W_j$ $(j\in I),$
take any $w\in W_j.$
By Hausdorffness of the coarse topology,
$x$ has a neighborhood $U$ in
that topology not containing $w,$ so
by~(\ref{d.basis+})(ii), it has such a neighborhood of the
form $U=\{x\}\cup\bigcup_{i\in J} W_i$ $(J\subseteq I_x).$
Since $U$ does not contain $w\in W_j,$ we have $j\notin J,$
so $U$ is disjoint from $W_j.$
With these modifications, the proof of Lemma~\ref{L.h_x},
and hence that of Theorem~\ref{T.main},
go over to the present hypotheses.
\end{proof}

To what spaces $X$ might one apply such a construction?
The set $\beta(\omega)$ has a natural compact Hausdorff
topology, under which it is the Stone-\v{C}ech compactification
of the discrete space $\omega.$
However, taking $X$ to be the former space and the
dense subset $Y$ to which we eventually restrict our functions
to be the latter does not seem a good candidate for our purposes.
For in the abovementioned compact topology on $\beta(\omega),$ the
neighborhoods of any $\U\in\beta(\omega),$ when intersected
with $\omega,$ give precisely the members of $\U.$
Since $\U$ is an ultrafilter,
if we tried  to ``puncture'' these neighborhoods further,
in a way that affected their intersections with $\omega,$
we would get some neighborhoods that intersected
$\omega$ in the empty set; i.e., $\omega$ would cease to be dense.

On the other hand, the less exotic space $2^c,$ i.e.,
the continuum power of the discrete space $2,$
also has countable dense subsets; let me sketch how to obtain one.
Identify $c$ with the set of $\!\{0,1\}\!$-valued functions on $\omega.$
Let $\r{Boole}(\omega)$ be the free Boolean algebra on
an $\!\omega\!$-tuple of indeterminates $x_0,\,x_1,\,\dots\,,$ i.e.,
the set of finitary Boolean operations in countably many variables.
Then to each $b(x_0,x_1,\dots)\in\r{Boole}(\omega)$
we can associate a subset of $c,$ namely the set $S_b$ of
$\!\omega\!$-tuples $(e_0,e_1,\dots)$ $(e_i\in\{0,1\})$
such that $b(e_0,e_1,\dots)=1.$
(Intuitively, the set of assignments
that ``satisfy'' the Boolean condition given by $b.)$
If we now think of $2^c$ as the power set of $c,$ so that
each $S_b$ is a member thereof, I claim that the countable set
$\{S_b\mid b\in\r{Boole}(\omega)\}\subseteq 2^c$ is dense.
Indeed, a basis of the topology of $2^c$ is given by the
solution-sets of statements saying that a certain finite list of
elements of $c$ should, and another finite list should
not, belong to the members of $2^c$ considered.
Since every finite family of elements of $c$ may be distinguished
by looking at finitely many coordinates, we can find
some $b\in\r{Boole}(\omega)$ such that $S_b$ satisfies the
criterion to belong to the given solution-set.

It should be possible to ``puncture'' the topology
on $2^c,$ so as to get stronger, non-compact topologies under which the
above countable set remains dense.
Whether one could get such a topology that
satisfied~(\ref{d.basis+}), or some other condition
from which one could prove that the group $G$ of
$\!\Z\!$-valued functions continuous in that topology had
dual free on the set of evaluations at
the points of $2^c,$ is not clear to me.

\section{Some other subgroups of $\Z^\omega.$}\label{S.examples}

We have seen that it is not likely that one can get
an affirmative answer to Question~\ref{Q.2^c} by finding
a group $G\subseteq\Z^\omega$ whose dual is spanned by evaluations
of elements at members of $\beta(\omega).$
However, while hoping to do so, I came upon some curious subgroups $G,$
which I sketch below for their own interest.

My idea was that since $\beta(\omega)$ is constructed
using only the set-theoretic structure of $\omega,$ and not
its order, etc., one should look at
groups whose definitions likewise ``treat
all points of $\omega$ alike''; i.e., subgroups of $\Z^\omega$
invariant under the action of the full symmetric group on $\omega.$
These contrast with subgroups of the sort commonly studied.
(For instance, those in \cite{Specker} consist of
the sequences with prescribed bounds on their growth rates.)

I will start with a class of examples that actually offers a faint
hope of giving a construction with the desired sort of dual.
Suppose we take any nondiscrete Hausdorff
\emph{group topology} $\mathcal{T}$
on $\Z$ (a topology under which the group operations are
continuous; for instance, the $\!p\!$-adic topology for some prime $p,$
or the topology induced by an embedding in the circle group,
or one of the topologies constructed in~\cite{IP+EZ}).
Now let $G_\mathcal{T}\subseteq\Z^\omega$ consist of all $g$ such
that the set of values of $g$ has compact closure $C(g)$ within
$\Z$ under that topology.
(For instance, if $(k_i)_{i\in\omega}$ is a sequence
of integers converging under $\mathcal{T}$ to an integer $k,$
then any $g\in\Z^\omega$ whose components all lie
in $\{k_i\mid i\in\omega\}\cup\{k\}$ will have this property.)
It is easy to see that $G_\mathcal{T}$ is a subgroup of $\Z^\omega.$
Given $g\in G_\mathcal{T},$ every ultrafilter $\U$
on $\omega$ induces an ultrafilter $g(\U)$
on $g(\omega)\subseteq C(g),$ which, by compactness of the latter
set, will converge to some element $\lim_{\U}g\in C(g)\subseteq\Z.$
For each $\U,$ this construction $g\mapsto\lim_{\,\U}\,g$ is a
homomorphism $h_\U:G_\mathcal{T}\to\Z.$
Whether there is a topology $\mathcal{T}$ on $\Z$ such that these
homomorphisms $h_\mathcal{U}$
span $\r{Hom}(G_\mathcal{T},\Z),$ I do not know.

A different sort of subgroup invariant under all permutations
of $\omega$ is
\begin{equation}\begin{minipage}[c]{35pc}\label{d.oEk}
$\{f\in\Z^\omega\mid(\exists\,k\in\omega)\,(\forall\,n>0)\ f$
takes on at most $k$ values modulo $n\}.$
\end{minipage}\end{equation}
A somewhat larger subgroup is
\begin{equation}\begin{minipage}[c]{35pc}\label{d.oEk_aa}
$\{f\in\Z^\omega\mid(\exists\,k\in\omega)\,(\forall\,n>0)\ f$ assumes
at most $k$ values modulo $n$ infinitely many times$\}.$
\end{minipage}\end{equation}
A much smaller subgroup of~(\ref{d.oEk_aa}) is
\begin{equation}\begin{minipage}[c]{35pc}\label{d.few_not_div}
$\{f\in\Z^\omega\mid(\forall\,n>0)$ there are only finitely
many $m\in\omega$ such that $n\not|f(m)\}.$
\end{minipage}\end{equation}

For a final class of examples, let us start with any
set $S$ of integers.
Then the subgroup $G_S\subseteq\Z^\omega$ generated by the set
$(\{0\}\cup S)^\omega$ is clearly
invariant under permutations of $\omega.$
It can be described as
\begin{equation}\begin{minipage}[c]{35pc}\label{d.bdd_nr_fr_S}
$G_S\ =\ \{f\in\Z^\omega\mid(\exists\,k\in\omega)$ every value assumed
by $f$ is the sum of at most $k$ terms taken from $\pm S$
(counting repetitions)$\}.$
\end{minipage}\end{equation}
So, for instance, if $S=\{1\},$ then $G_S$ is
the group $B$ of bounded functions; if $S$ is the set of powers
of $2,$ then $G_S$ can be described as the set of all sequences
of integers whose binary expressions
(ignoring initial $\pm$ signs) have a common bound on the number
of substrings ``$10$'' that they contain.
If $S$ has compact closure under some group
topology $\mathcal{T}$ on $\Z,$
then $G_S$ will be contained in the
group $G_{\mathcal T}$ discussed above.
On the other hand, if $S$ is not sufficiently sparse, e.g., if it is
the set of all squares, then $G_S$ is the full group~$\Z^\omega.$

\end{document}